\documentstyle[11pt]{article}
\begin{document}
\setlength{\oddsidemargin}{0cm}
\setlength{\evensidemargin}{0cm}
\baselineskip=20pt

\begin{center} {\Large\bf
Left-symmetric Algebras From Linear Functions}  \end{center}

\bigskip

\begin{center}  { \large Chengming ${\rm Bai}^{1,2,3}$} \end{center}

\begin{center}{\it 1. Nankai Institute of
Mathematics, Tianjin 300071, P.R. China } \end{center}

\begin{center}{\it 2. Liu Hui Center for Applied Mathematics, Tianjin
300071, P.R. China}\end{center}

\begin{center}{\it 3. Dept. of Mathematics, Rutgers, The State University of
New Jersey, Piscataway, NJ 08854, U.S.A.}\end{center}

\vspace{0.3cm}

\begin{center} {\large\bf   Abstract } \end{center}

In this paper, some left-symmetric algebras are constructed from linear
functions. They
include a kind of simple left-symmetric algebras and some examples
appearing in mathematical physics. Their complete classification is also
given, which shows that they can be regarded as generalization of
certain 2-dimensional left-symmetric algebras.

\vspace {0.2cm}


\vspace{0.2cm}


\bigskip

{\bf I. Introduction}

\vspace{0.3cm}

A left-symmetric algebra is an algebra whose associator is left-symmetric:
Let $A$ be a vector space over a field ${\bf F}$ with a bilinear product
$(x,y)\rightarrow xy$. $A$ is called a left-symmetric algebra if for any
$x,y,z\in A$, the associator
$$(x,y,z)=(xy)z-x(yz)\eqno (1.1)$$
is symmetric in $x,y$, that is,
$$(x,y,z)=(y,x,z),\;\;{\rm or}\;\;{\rm
equivalently}\;\;(xy)z-x(yz)=(yx)z-y(xz).\eqno (1.2)$$

Left-symmetric algebras are a class of non-associative algebras arising from
the study of convex homogenous cones, affine manifolds and affine
structures on Lie groups([V], [Ki], [P]). Moreover, they have very close
relations with
many problems in mathematical physics. For example, they appear as an
underlying structure of those Lie algebras that possess a phase space
([Ku1-4], thus they form a natural category from the point of view of
classical and quantum mechanics) and there is a close relation between them
and classical Yang-Baxter equation([ES],[ESS],[GS]).

However, due to the non-associativity, there is not a suitable
representation theory of left-symmetric algebras. It is also known that the
definition identity (1.2) of left-symmetric algebras involves the quadric
forms of structure constants, which is not linear in general ([Ki]). Hence
it is quite difficult to study them. Therefore one of the most important
problems is how to construct interesting left-symmetric algebras. One way is
to construct them through some well-known algebras and algebraic structures.
This can be regarded as
a kind of ``realization theory''. For example, there is a study of
realization of Novikov algebras (they are left-symmetric algebras with
commuting right multiplications) from commutative associative algebras and
Lie algebras in Refs. [BM3-5].
Another way is to try to reduce the ``non-linearity'' in certain sense.
Combining these two ways, a natural and simple way is to construct
left-symmetric algebras from linear
functions, which is the main content of this paper.

On the other hand, there are many examples of left-symmetric algebras
appearing in mathematical physics ([BN],[GD],[GS],[SS],etc.). For example,
let $V$ be
a vector space over the complex field ${\bf C}$ with the ordinary scalar
product $(,)$ and $a$
be a fixed vector in $V$, then
$$u*v=(u,v)a+(u,a)v,\forall u,v\in V,\eqno (1.3)$$
defines a left-symmetric algebra on $V$ which gives the integrable
(generalized) Burgers equation([SS], [S])
$$U_t=U_{xx}+2U*U_x+(U*(U*U))-((U*U)*U).\eqno (1.4)$$
However, such examples are often scattered and independent in different
references of mathematical physics. And in most of the cases, there is
neither a good mathematical motivation nor a further study. In this paper,
our construction not only has a natural motivation from the point of view of
mathematics, but also can be regarded as a kind of generalization of the
examples given by equation (1.3). Moreover, a systematic study is given.

The algebras that we consider in this paper are of finite dimension and over
${\bf C}$. The paper is organized as follows. In section (II), we construct
left-symmetric algebras from linear functions. In section
(III), we give their classification. In section (IV), we
discuss some properties of these left-symmetric algebras and certain
application in mathematical physics.

\newpage

{\bf II. Constructing left-symmetric algebras from linear functions}

\vspace{0.3cm}

Let $A$ be a vector space in dimension $n$. In general, we assume $n\geq 2$.
Just as said in the
introduction, motivated by the study of algebraic structure itself and some
equations in integrable systems, it is natural to consider the
left-symmetric algebras
satisfying the following conditions: for any two vectors $x,y$ in $A$, the
product $x*y$ is still in the subspace spanned by $x,y$, that is, any two
vectors make up a subalgebra in $A$. Thus, it is natural to assume
$$x*y=f_1(x,y)x+f_2(x,y)y,\;\;\forall x,y\in A,\eqno (2.1)$$
where $f_1, f_2:A\times A\rightarrow {\bf C}$ are two functions. In general,
$f_1$ and $f_2$ are not necessarily linear. However, if they are not linear
functions, they cannot be decided by their values at a basis of $A$. Hence
the problem turns to be more complicated, even more complicated than the
study of the algebra itself.

Therefore, we can assume that $f_1$ and $f_2$ are linear functions. Since
the
algebra product $*$ is bilinear, for $f_1\ne 0$, $f_1$ depends on only
$y$, that is, $f_1$ is not a linear function depending on $x$. Otherwise,
for any $\lambda \in {\bf C}$, we have
$$(\lambda x)*y=f_1(\lambda x, y)\lambda x+f_2(\lambda
x,y)y=\lambda^2f_1(x,y)x+\lambda f_2(x,y)y=\lambda
(f_1(x,y)x+f_2(x,y)y).\eqno (2.2)$$
Hence $f_1(x,y)=0,\;\forall x,y\in A$, which is a contradiction. Similarly,
$f_2$ depends on only $x$. Thus, we can set $f_1(x,y)=f(y),f_2(x,y)=g(x)$,
where $f,g:A\rightarrow {\bf C}$ are two linear functions.

{\bf Proposition 2.1}\quad Let $A$ be a vector space in dimension $n\geq 2$.
Let
$f,g:A\rightarrow {\bf
C}$ be two linear functions. Then the product
$$x*y=f(y)x+g(x)y,\;\;\forall x,y\in A\eqno (2.3)$$
defines a left-symmetric algebra if and only if $f=0$ or $g=0$. Moreover,
when $f=0$ or $g=0$, the above equation defines an associative algebra.

{\bf Proof}\quad  For any $x,y,z\in A$, the associator
\begin{eqnarray*}
(x,y,z)&=&(x*y)*z-x*(y*z)\\
&=&f(y)(f(z)x+g(x)z)+g(x)(f(z)y+g(y)z)-f(z)(f(y)x+g(x)y)\\
&&-g(y)(f(z)x+g(x)z)\\
&=&f(y)g(x)z-g(y)f(z)x
\end{eqnarray*}
Hence $(x,y,z)=(y,x,z)$ if and only if for any $y,z\in A$, $g(y)f(z)=0$,
that is, $f=0$ or
$g=0$. Moreover, when $f=0$ or $g=0$, $(x,y,z)=0$. Thus the proposition
holds.
\hfill $\Box$

Let $L_x,\;R_x$ denote the left and right multiplication respectively,
i.e., $L_x(y)=xy,\;R_x(y)=yx,\;\forall x,\;y\in A$.

{\bf Corollary 2.2}\quad With the conditions in above proposition, we have

(1) If $f=0, g\ne 0$, then there exists a basis $\{e_1,\cdots, e_n\}$ in $A$
such that $L_{e_1}={\rm Id}, L_{e_i}=0,i=2,\cdots,n$, where ${\rm Id}$ is
the identity transformation.

(2) If $g=0, f\ne 0$, then there exists a basis $\{e_1,\cdots, e_n\}$ in $A$
such that $R_{e_1}={\rm Id}, R_{e_i}=0,i=2,\cdots,n$.

(3) If $f=g=0$, then $A$ is a trivial algebra, that is, all products are
zero.

{\bf Proof}\quad For any linear function $g:A\rightarrow {\bf C}$, if $g\ne
0$, due to the linearity of $g$ and the direct sum of vector spaces
$$A={\rm Ker} g\oplus g(A)={\rm Ker} g\oplus {\bf C},$$
there exists a basis  $\{e_1,\cdots, e_n\}$ in $A$ such that $g(e_1)\ne 0,
g(e_i)=0, i=2,\cdots,n$. Furthermore, we can normalize $g$ by $g(e_1)=1$.
Hence (1) and (2) follows. (3) is obvious. \hfill $\Box$

{\bf Remark 1}\quad There is a natural matrix representation of above
associative algebras ([Bu2]). Let $\{E_{ij}\}$ be the canonical basis of
$gl(n)$, that is, $E_{ij}$ is a $n\times n$ matrix with 1 at $i$th row and
$j$th column and zero at other places. Then the algebra in above case (1)
(respectively (2)) is an associative subalgebra of $gl(n)$ (under the
ordinary matrix product) with $e_i=E_{1i}$ (respectively $e_i=E_{i1}$).

It is well known that the commutator of a left-symmetric algebra
$[x,\;y]=xy-yx$ defines a (sub-adjacent) Lie algebra ([Ki], [P], etc.).

{\bf Corollary 2.3}\quad The sub-adjacent Lie algebras of the associative
algebras defined by equation (2.3) with $g=0,f\ne 0,$ or $f=0, g\ne 0$ are
isomorphic
to the following 2-step solvable Lie algebra:
$$A=<e_i, i=1,\cdots,n|[e_1,e_i]=e_i,\;i=2,\cdots, n,\;{\rm other}\;\;{\rm
products}\;\;{\rm
are}\;\;{\rm zero}>.\eqno (2.4)$$

{\bf Proof}\quad For case (1) in Corollary 2.2, the conclusion is obvious.
For case (2) in Corollary (2.2), we only need a linear
transformation by letting $e_1$ be $-e_1$ and $e_i$ still be $e_i$
($i=2,\cdots,n$), which the conclusion follows.
\hfill $\Box$

{\bf Remark 2}\quad The above conclusion also can be obtained from equation
(2.3) directly. That is, the Lie algebra given by
$[x,y]=(f-g)(x)y-(f-g)(y)x$ is isomorphic to the Lie algebra given by
equation (2.4) for $g\ne f$. In fact, this algebra can be regarded as a
(unique!) non-abelian Lie algebra constructed from linear functions: it is
easy to show
that the product $[x,y]=f(x)y+g(y)x$ defines a Lie algebra if and only if
$f(x)=-g(x),\forall x\in A$.

Due to the above discussion, in order to get non-associative left-symmetric
algebras, we need to extend the above construction. A simple extension of
equation
(2.3) is to add a fixed vector $c\ne 0$ as follows:
$$x*y=f(x)y+g(y)x+h(x,y)c,\;\;\forall x,y\in A,\eqno (2.5)$$
where $h:A\times A\rightarrow {\bf C}$ is a non-zero bilinear function. The
above equation
(2.5) can be understood that for any two vectors $x,y$, the three vectors
$x,y,c$ make up a
subalgebra in $A$. Moreover, if $h$ is symmetric, then its sub-adjacent Lie
algebra
is isomorphic to the Lie algebra given by equation (2.4) ($f\ne g)$ or the
abelian Lie algebra ($f=g$).

For a further study, we give a lemma on linear functions at first:

{\bf Lemma 2.4}\quad Let $A$ be a vector space in dimension $n\geq 2$.
Let $f,g:A\rightarrow {\bf C}$ be two linear functions and $h:A\times
A\rightarrow {\bf C}$ be a symmetric bilinear function.

(1) If for any $x,y\in A$,
$f(x)g(y)=f(y)g(x)$, then $f=0$, or $g=0$, or $f\ne 0, g\ne 0$ and there
exists $\alpha\in {\bf C}, \alpha\ne 0$ such that $f(x)=\alpha g(x),\forall
x\in A$.

(2) If for
any $x,y,z\in A$, $f(x)h(y,z)=f(y)h(x,z)$, then $f=0$, or $h=0$, or there
exists a basis $\{ e_1,\cdots, e_n\}$ in $A$ and $\alpha\in {\bf C},
\alpha\ne 0$ such that $f(x)=h(x,\alpha e_1),\forall x\in A; h(e_1,e_1)= 1,
h(e_i, e_j)=0$, $i=2,\cdots, n, j=1,\cdots, n$.

{\bf Proof}\quad For a linear function $f$, if $f\ne 0$, we can choose a
basis $\{ e_1,\cdots, e_n\}$ in $A$ such that $f(e_1)\ne 0,
f(e_2)=\cdots=f(e_n)=0$. If $g\ne 0$, then from $f(x)g(y)=f(y)g(x)$, we can
have $g(e_1)\ne 0, g(e_i)=0,i=2,\cdots,n$.
Let $\alpha=f(e_1)/g(e_1)$, then by linearity, for any $x\in A$, we have
$f(x)=
\alpha g(x)$.

Similarly, for $f\ne 0$ and the basis $\{ e_1,\cdots, e_n\}$ in $A$ such
that $f(e_1)\ne 0,
f(e_2)=\cdots=f(e_n)=0$, we have $h=0$ or $h(e_1,e_1)\ne 0$ and $h(e_i,
e_j)=0, i=2,\cdots, n, j=1,\cdots, n$. For the latter case, we can normalize
$h$ by
$h(e_1,e_1)=1$. Thus, we still have $f(x)=h(x, \alpha e_1),\forall x\in A$,
where $\alpha=\frac{f(e_1)}{h(e_1,e_1)}=f(e_1)$. \hfill $\Box$

{\bf Theorem 2.5}\quad With the conditions in above lemma and $h\ne 0$,
equation (2.5) defines a left-symmetric algebra if and only if the
functions $f,g,h$ belong to one of the following cases:

(1) $f=g=0$, $h(x,c)=0,\forall x\in A$;

(2) $f=g=0$, and there exists a basis $\{ e_1,\cdots,e_n\}$ such that
$h(e_1,e_1)= 1, h(e_i, e_j)=0$, $i=2,\cdots, n,
j=1,\cdots, n$ and $c=\sum\limits_{i=1}^na_ie_i$ with $a_1\ne 0$;

(3) $g=0, f\ne 0$, $f(x)=h(x,c),\;\forall x\in A$;

(4) $g=0$, $f\ne 0$, and there exists a basis $\{e_1,\cdots,
e_n\}$ and $\alpha\in {\bf C}, \alpha\ne 0$ such that $f(x)=h(x,c-\alpha
e_1), h(e_1,e_1)=1, h(e_i,e_j)=0,i=2,\cdots,n,j=1,\cdots,n$ and
$c=\sum\limits_{i=1}^na_ie_i$ with $a_1\ne \alpha$;

(5) $f=0, g\ne 0$, $g(x)=-h(x,c), \forall x\in A$ and $h(c,c)=0$;

(6) $f=0, g\ne 0$, $h(x,c)=0,\forall x\in A$, and there exists a basis
$\{e_1,\cdots, e_n\}$ and $\alpha\in {\bf C}, \alpha\ne 0$ such that
$g(x)=h(x,\alpha e_1), h(e_1,e_1)=1,
h(e_i,e_j)=0,i=2,\cdots,n,j=1,\cdots,n$;

(7) $f\ne 0, g\ne 0$, $f(c)\ne 0$ and there exists $\alpha\in {\bf C},
\alpha\ne
0$ such that $g(x)=\alpha f(x), h(x,y)=-\frac{f(x)f(y)}{f(c)},\forall x\in
A$.

{\bf Proof}\quad For any $x,y,z\in A$, the associator
\begin{eqnarray*}
(x,y,z)&=&(x*y)*z-x*(y*z)\\
&=&
f(x)[f(y)z+g(z)y+h(y,z)c]+g(y)[f(x)z+g(z)x+h(x,z)c]\\
&&+h(x,y)[f(c)z+g(z)c+h(c,z)c]-f(y)[f(x)z+g(z)x+h(x,z)c]\\
&&-g(z)[f(x)y+g(y)x+h(x,y)c]-h(y,z)[f(x)c+g(c)x+h(x,c)c]\\
&=&
[-f(y)g(z)-g(c)h(y,z)]x+[g(y)f(x)+f(c)h(x,y)]z\\
&&+[g(y)h(x,z)-f(y)h(x,z)+h(x,y)h(c,z)-h(y,z)h(x,c)]c
\end{eqnarray*}
Then by left-symmetry, we can get the following equations: for any $x,y,z\in
A$:
$$f(y)g(x)=g(y)f(x);\eqno (2.6)$$
$$f(y)g(z)+g(c)h(y,z)=0;\eqno (2.7)$$
$$[(g-f)(y)+h(y,c)]h(x,z)=[(g-f)(x)+h(x,c)]h(y,z).\eqno (2.8)$$
From equation (2.6) and using Lemma 2.4, we can consider the following
cases:

Case (I): $f=g=0$. There is only one non-trivial equation
$h(y,c)h(x,z)=h(x,c)h(y,z)$. Let $h'(x)=h(x,c)$, then by Lemma 2.4, we know
that $h'(x)=0$ or there
exists a basis $\{ e_1,\cdots, e_n\}$ in $A$ and $\alpha\in {\bf C},
\alpha\ne 0$ such that $h'(x)=h(x,\alpha e_1),\forall x\in A; h(e_1,e_1)= 1,
h(e_i, e_j)=0$, $i=2,\cdots, n, j=1,\cdots, n$. The former is the case (1)
and the latter is the case (2) since $h(x,c)=h(x,\alpha e_1)$ implies that
$a_1=\alpha\ne 0$ for $c=\sum\limits_{i=1}^na_ie_i$.

Case (II): $g=0, f\ne0$. Then equation (2.7) is satisfied.  From equation
(2.8) and using lemma 2.4 again, we have $f(x)=h(x,c)$ or there exists a
basis $\{e_1,\cdots, e_n\}$ such that $h(e_1,c)-f(e_1)\ne 0,
f(e_i)=h(e_i,c)=0;\;h(e_1,e_1)= 1, h(e_i,e_j)=0,i=2,\cdots,n,j=1,\cdots,n$.
The former is the case (3) and the latter is the case (4) where
$\alpha=-f(e_1)+h(e_1,c)$. Notice for the latter, $f\ne 0$ if and only if
$a_1\ne \alpha$ for $c=\sum\limits_{i=1}^na_ie_i$.

Case (III): $f=0, g\ne 0$. From
equation (2.7), we have $g(c)=0$. As the same as the discussion in
Case (II), equation (2.8) implies that $g(x)=-h(x,c)$ or there exists a
basis $\{e_1,\cdots, e_n\}$ such that $g(e_1)+h(e_1,c)\ne 0,
g(e_i)=h(e_i,c)=0,\;h(e_1,e_1)= 1, h(e_i,e_j)=0,i=2,\cdots,n,j=1,\cdots,n$.
The former is the case (5). For the latter, we have $g(x)=-h(x,c)+\alpha
h(x,e_1)$ where $\alpha=g(e_1)+h(e_1,c)$. Set $c=\sum\limits_{i=1}^n
a_ie_i$, then $g(c)=-a_1^2+\alpha a_1=0$. Thus $a_1=\alpha$ or $a_1=0$.
Therefore if $g\ne 0$, we have $h(x,c)=0$ and $g(x)=h(x,\alpha e_1)$ which
is just the case (6).

Case (IV): $f\ne 0$, $g\ne 0$. Thus there exists $\alpha \ne 0$ such
that $g(x)=\alpha f(x)$. Hence from equation (2.7) and the assumption
$h\ne 0$, we know that $f(c)\ne 0$ and
$h(x,y)=-\frac{f(x)f(y)}{f(c)},\forall x\in A$. It is easy to see that
equation (2.8) holds under these conditions. This is the case (7). \hfill
$\Box$

{\bf Corollary 2.6}\quad The left-symmetric algebras given in Theorem 2.5
are commutative (hence associative), if and only if their sub-adjacent Lie
algebras are abelian, if and only if they belong to the case (1), (2) and
(7)
with $\alpha=1$.

By direct checking, we have

{\bf Corollary 2.7}\quad Let $A$ be a left-symmetric algebra in Theorem 2.5.

(1) If $A$ is in the case (1), (2), (4), (6), (7), then the corresponding
bilinear function $h$ satisfies
$$h(x*y,z)=h(y*x,z)=h(x*z, y),\;\;\forall x,y,z\in A.\eqno (2.9)$$

(2) If $A$ is in the case (3), then the corresponding bilinear function $h$
is invariant under the
following sense:
$$h(x*y,z)=h(x,z*y),\forall x,y,z\in A.\eqno (2.10)$$
That is, for every $x,y,z\in A$, $h(R_x(y),z)=h(y,R_x(z))$ ($R_x$ is
self-adjoint).

(3) If $A$ is in the case (5), then the corresponding bilinear function $h$
satisfies
$$h(x*y,z)+h(y,x*z)=0,\forall x,y,z\in A.\eqno (2.11)$$
That is, for every $x,y,z\in A$, $h(L_x,z)+h(x,L_z)=0$.

\newpage

{\bf III.\quad Classification of left-symmetric algebras from linear
functions}

\vspace{0.3cm}

In this section, we discuss the classification of left-symmetric algebras
given in Theorem 2.5. Since the bilinear function $h$ appearing in the case
(2),
(4), (6) and (7) has been (almost) decided completely, we give the
classification of these cases at first.

{\bf Proposition 3.1}\quad Let $A$ be a left-symmetric algebra in the case
(2)
with dimension $n\geq 2$. Then $A$ is isomorphic to the following
algebra:(we only give the non-zero products)
$$A_{(2)}=< e_i,i=1,\cdots,n|e_1e_1=e_1>.\eqno (3.1)$$

{\bf Proof}\quad For $c=\sum\limits_{i=1}^n a_ie_i$ with $a_1\ne 0$, let
$$e_1'=\frac{1}{a_1}e_1+\frac{1}{a_1^2}\sum\limits_{i=2}^n a_ie_i; \;\;
e_j'=e_j,j=2,\cdots,n,$$
then under the new basis, equation (3.1) follows. \hfill $\Box$

{\bf Proposition 3.2}\quad Let $A$ be a left-symmetric algebra in the case
(4)
with dimension $n\geq 2$. Then $A$ is isomorphic to one of the following
algebras:
$$A^1_{(4)}=< e_i,i=1,\cdots,n|e_1e_1=e_1+e_2, e_1e_j=e_j,
j=2,\cdots,n>;\eqno
(3.2)$$
$$A^\lambda_{(4)}=< e_i,i=1,\cdots,n|e_1e_1=\lambda e_1, e_1e_j=e_j,
j=2,\cdots,n;>,\;\;\lambda\ne
1,2.\eqno (3.3)$$

{\bf Proof}\quad For the case (4), we have
$$e_1*e_1=h(e_1,(a_1-\alpha)e_1)e_1+h(e_1,e_1)c=(a_1-\alpha)e_1+c=(2a_1-\alpha)e_1+\sum_{i=2}^n
a_ie_i,$$
$$e_1*e_i=(a_1-\alpha)e_i,\;\;e_i*e_j=0,i=2,\cdots,n,j=1,\cdots,n.$$
If $a_1=0$, then $c=\sum\limits_{i=2}^n a_ie_i\ne 0$. Without losing
generality, we suppose $a_2\ne 0$.
Let
$$e_1'=-\frac{1}{\alpha}e_1;\;\; e_2'=\frac{1}{\alpha^2}c; \;\;
e_j'=e_j,\;j=3,\cdots,n,$$
then under the new basis, we can get equation (3.2).

\noindent If $a_1\ne 0$ and $a_1\ne \alpha$, then let
$$e_1'=\frac{1}{a_1-\alpha}e_1+\frac{1}{(a_1-\alpha)a_1}\sum_{i=2}^n
a_ie_i;\;\;e_j'=e_j,j=2,\cdots,n.$$
Hence under the new basis, we have
$$e_1'*e_1'=\frac{2a_1-\alpha}{a_1-\alpha}e_1',
e_1'*e_i'=e_i',\;\;e_i'*e_j'=0,i=2,\cdots,n,j=1,\cdots,n.$$
Set $\lambda=\frac{2a_1-\alpha}{a_1-\alpha}$ which gives equation (3.3).
Notice that $\lambda \ne 1, 2$ since $a_1\ne 0, a_1\ne \alpha$ and
$\alpha\ne 0$.

\mbox{}\hfill $\Box$

As the same as the proof of equation (3.2), we have

{\bf Proposition 3.3}\quad Let $A$ be a left-symmetric algebra in the case
(6)
with dimension $n\geq 2$. Then $A$ is isomorphic to
$$A_{(6)}=< e_i,i=1,\cdots,n|e_1e_1=e_1+e_2, e_je_1=e_j, j=2,\cdots,n>.\eqno
(3.4)$$

{\bf Proposition 3.4}\quad Let $A$ be a left-symmetric algebra in the case
(7)
with dimension $n\geq 2$. Then $A$ is isomorphic to one of the following
algebras:
$$A^{\alpha}_{(7)}=< e_i,i=1,\cdots,n|e_1e_1=\alpha e_1, e_1e_j=e_j,
e_je_1=\alpha e_j,
j=2,\cdots,n>,\;\;\alpha\ne 0.\eqno (3.5)$$

{\bf Proof}\quad Without losing generality, we can choose a basis $\{
e_1,\cdots,e_n\}$ such that $e_1=c$ and $f(e_2)=\cdots=f(e_n)=0$. Hence
$$e_1*e_1=\alpha f(e_1)e_1,\;e_1*e_j=f(e_1)e_j,\;e_j*e_1=\alpha
f(e_1)e_j,\;j=2,\cdots,n.$$
The conclusion follows by the basis transformation
$e_1'=\frac{1}{f(e_1)}e_1, e_j'=e_j,j=2,\cdots,n$.\hfill $\Box$

In order to classify the left-symmetric algebras in other cases, we
need the following lemma:

{\bf Lemma 3.5}\quad Let $A$ be a finite-dimensional algebra over
${\bf C}$. Let $A=A_1\oplus A_2$ as the direct sum of two
subspaces and $A_1$ be a subalgebra.  Assume that, for every $x\in
A_1$, $L_x$ and $R_x$ acts on $A_2$ is zero or ${\rm Id}$. If
there exists a non-zero vector $v\in A_1$ such that for any two
vectors $x,y\in A_2$, $xy=yx\in {\bf C}v$, then the classification
of the algebraic operation in $A_2$ (without changing other
products) is given by the classification of symmetric bilinear
forms on a $n$-dimensional vector space over ${\bf C}$, where
$n=\dim A_2$. That is, there exists a basis $\{ e_1, \cdots,
e_n\}$ in $A_2$ such that the classification is given as follows:
$A_2$ is trivial or for every $k=1,\cdots, n$:
$$e_ie_j=\delta_{ij}v,\;\; i,j=1,\cdots,k;\;\;e_ie_j=0,\;\;{\rm
otherwise}.\eqno
(3.6)$$

{\bf Proof}\quad From the assumption, there exists a symmetric bilinear form
$f:A_2\times A_2\rightarrow {\bf C}$ such that
$$xy=f(x,y)v,\;\;\forall x,y\in A_2.$$
Moreover, any linear transformation of $A_2$ does not change the operation
relations between $A_1$ and $A_2$, hence the whole algebra $A=A_1\oplus
A_2$.
Thus the classification of $A_2$ is decided completely by the classification
of symmetric bilinear forms on a vector space in dimension $\dim A_2$.
Therefore
there exists a basis $\{ e_1, \cdots, e_n\}$ in $A_2$ such that the matrix
$(f(e_i,e_j))$ is zero or a diagonal matrix with the first $k$
($k=1,\cdots,n$)
elements are 1 and the others are zero on the diagonal, which gives equation
(3.6). It is easy to show that for different $k$, the algebras are not
mutually isomorphic.
\hfill $\Box$

{\bf Proposition 3.6}\quad The classification of left-symmetric algebras in
the
case (1) with dimension $n\geq 2$ is equivalent to the classification of
symmetric bilinear forms on a $(n-1)$-dimensional vector space. The
classification is given as follows: for every $k=0,\cdots, n-1$,
$$A_{(1)}^{(k)}=< e_i,i=1,\cdots,n|e_je_j=e_1, j=2,\cdots,
k+1>.\eqno (3.7)$$

{\bf Proof}\quad Let $A$ be a left-symmetric algebra in the case (1) with
dimension $n\geq 2$. We can choose a basis $\{e_1,\cdots, e_n\}$ such that
$e_1=c$. Thus we have
$$e_1*e_1=e_1*e_j=e_j*e_1=0, \;\;e_j*e_k=h(e_j,e_k)e_1,\;j,k=2,\cdots, n.$$
Let $A_1$ be a subspace spanned by $e_1$ and $A_2$ be a subspace spanned by
$e_2,\cdots, e_n$. Then by Lemma 3.5, the proposition holds. \hfill $\Box$

{\bf Proposition 3.7}\quad The classification of left-symmetric algebras in
the
case (3) with dimension $n\geq 2$ is given by the following matrices
($F=(h(e_i,e_j))$,
where $\{e_1, \cdots,e_n\}$ is a basis)
$$F_1=I;\;\;
F_2^{(k)}=\left( \matrix{1 & 0 &0\cr 0 &0 &0 \cr 0 &0&
A^{(k)}\cr}\right);\;\;
F_3^{(k)}=\left( \matrix{0 & 1 &0\cr 1 &0 &0 \cr 0 &0& A^{(k)}\cr}\right)\;
,\eqno (3.8)$$
where $A^{(k)}={\rm diag} (1,\cdots 1,0,\cdots 0)$ is a $(n-2)\times
(n-2)$ diagonal matrix with the first $k$ elements are 1 and the others are
zero on the diagonal, $k=0,1\cdots, n-2$. The corresponding left-symmetric
algebras
are
$$A_{(3)}^1=< e_i,i=1,\cdots,n|e_1e_1=2e_1, e_1e_j=e_j,
e_je_j=e_1,j=2,\cdots,n>;\eqno (3.9)$$
$$A_{(3)}^{(k),2}=< e_i,i=1,\cdots,n|e_1e_1=2e_1, e_1e_j=e_j, e_le_l=e_1,
j=2,\cdots,n, l=3,\cdots, k+2>;\eqno (3.10)$$
$$A_{(3)}^{(k),3}=< e_i,i=1,\cdots,n|e_1e_2=e_1, e_2e_1=2e_1, e_2e_2=e_2,
e_2e_j=e_j, e_le_l=e_1, $$
$$\mbox{}\hspace{5cm} j=3,\cdots,n, l=3,\cdots, k+2>.\eqno (3.11)$$

{\bf Proof}\quad Let $A$ be a left-symmetric algebra in the case (3) with
dimension $n\geq 2$. Without losing generality, we can assume $c=e_1$. At
first we consider the case $h(e_1,e_1)\ne 0$. Thus we can choose
$e_2,\cdots, e_n$ such that $\{e_1,\cdots,e_n\}$ is a basis and
$h(e_1,e_j)=0,j=2,\cdots,n$. Set $h_{ij}=h(e_i,e_j)$. Therefore the product
of $A$ is given by
$$e_1*e_1=2h_{11}e_1,\;\; e_1*e_j=h_{11}e_j,\;\; e_j*e_1=0,\;\;
e_j*e_l=h_{jl}e_1,\;
j,l=2,\cdots,n.$$
Moreover, we can assume $h_{11}=1$ by letting $e_1'=\frac{1}{h_{11}}e_1,
e_j'=\frac{1}{\sqrt{h_{11}}}e_j,j=2,\cdots,n$.
Let $A_1={\bf C}e_1$ and $A_2$ be a subspace spanned by $e_2,\cdots, e_n$,
then from
Lemma 3.5, we know the classification of above algebras is just given by the
matrix $F_1$ and $F_2^{(k)}$ respectively, which corresponds to the
left-symmetric algebra
given by equation (3.9) and equation (3.10) respectively.

Next assume $h(c,c)=(e_1,e_1)=0$. Since there exists an element $u\in A$
such that
$h(u,c)\ne 0$, we can let $u=e_2$. Then we can choose $e_3,\cdots, e_n$ such
that $\{e_1,\cdots,e_n\}$ is a basis and $h(e_1,e_j)=0,j=1,3,\cdots,n$.
Hence we have
$$e_1*e_1=0, e_1*e_2=h_{12}e_1, e_2*e_1=2h_{12}e_1,
e_2*e_2=h_{12}e_2+h_{22}e_1,$$
$$ e_j*e_1=e_1*e_j=0, e_2*e_j=h_{12}e_j+h_{2j}e_1, e_j*e_2=h_{2j}e_1,
e_j*e_l=h_{jl}e_1, j,l=3,\cdots,n.$$
Let
$$e_1'=e_1;\;\; e_2'=\frac{1}{h_{12}}e_2-\frac{h_{22}}{2h_{12}}e_1;\;\;
e_j'=e_j-\frac{h_{2j}}{h_{12}^2}e_1,\;j=3,\cdots,n.$$
Under the new basis, we have
$$e_1'*e_1'=0, e_1'*e_2'=e_1', e_2'*e_1'=2e_1', e_2'*e_2=e_2',$$
$$ e_j'*e_1'=e_1'*e_j'=0, e_2'*e_j'=e_j', e_j'*e_2'=0, e_j'*e_l'=h_{jl}e_1',
j,l=3,\cdots,n.$$
Let $A_1$ be a subspace spanned by $e_1,e_2$ and $A_2$ be a subspace spanned
by $e_3,\cdots, e_n$,
then from Lemma 3.5, we know the classification of above algebras is just
given by the matrix $F_3^{(k)}$, which corresponds to the left-symmetric
algebra given by equation (3.11).\hfill $\Box$

As the same as the proof of the case $A_{(3)}^{(k),3}$ in above proposition,
we have

{\bf Proposition 3.8}\quad The classification of left-symmetric algebras in
the
case (5) with dimension $n\geq 2$ is given by the matrix $F_3^{(k)}$.
The corresponding left-symmetric algebras is ($k=0,1\cdots, n-2$)
$$A_{(5)}^{(k)}=< e_i,i=1,\cdots,n|e_2e_1=-e_1, e_2e_2=e_2, e_je_2=e_j,
e_le_l=e_1, 3\leq j\leq n, 3\leq l \leq k+2>.\eqno (3.12)$$

{\bf Corollary 3.9}\quad Let $A$ be a left-symmetric algebra in dimension
$n\geq 2$ given in Theorem 2.5. If the bilinear function $h$ is
non-degenerate, then $A$ is isomorphic to one of the following algebras:
$A_{(3)}^1$; $A_{(3)}^{(n-2),3}$; $A_{(5)}^{(n-2)}$.

{\bf Theorem 3.10}\quad When the dimension $n=2$, the left-symmetric
algebras
given in Theorem 2.5 are not (mutually) isomorphic except for
$$A_{(3)}^{(0),3}\sim A_{(7)}^{\frac{1}{2}};\;\;A_{(5)}^{(0)}\sim
A_{(4)}^{-1}.\eqno (3.13)$$
Moreover, with the associative algebras given in Corollary 2.2 together,
they include all 2-dimensional non-commutative left-symmetric algebras.

{\bf Proof}\quad Comparing the classification of 2-dimensional
left-symmetric algebras given in [BM1] or [Bu2], the conclusion follows
immediately. Notice that $A_{(3)}^{(0),3}$ is isomorphic to
$A_{(7)}^{\frac{1}{2}}$ by $e_1\rightarrow e_2, e_2\rightarrow 2e_1$ and
$A_{(5)}^{(0)}$ is isomorphic to  $A_{(4)}^{-1}$ by $e_1\rightarrow e_2,
e_2\rightarrow -e_1$, which the order of $e_1,e_2$ is changed
respectively.\hfill $\Box$

{\bf Remark 3}\quad Obviously, some commutative associative algebras such as
the direct sum of two fields ${\bf C}\oplus {\bf C}$ are not included in
above algebras. Moreover, we would like to point out that the above
conclusion is not obvious since for a general algebra, the ``linear''
construction like in this paper has certain restriction conditions for the
corresponding structure constants, which could not contain all (non-trivial)
examples.

{\bf Corollary 3.11}\quad When $n>2$, the left-symmetric algebras given in
Theorem 2.5 and Corollary 2.2 are not mutually isomorphic.

{\bf Proof}\quad It is easy to see that when $n>2$,
$A_{(3)}^{(k), 3}$ is not isomorphic to $ A_{(7)}^{\frac{1}{2}}$ and
$A_{(5)}^{(k)}$ is not isomorphic to  $A_{(4)}^{-1}$. With the special roles
of
$e_1,e_2$ in the algebraic operation and similar to the classification of
2-dimensional left-symmetric algebras in [BM1] or [Bu2], the conclusion
follows by a straightforward analysis. \hfill $\Box$

\vspace{0.3cm}

{\bf IV\quad Further Discussion}

\vspace{0.3cm}

In this section, we discuss some properties of the algebras given in the
previous sections and certain application in mathematical application.

{\bf Theorem 4.1}\quad The left-symmetric algebras given by equation (1.3)
are isomorphic to the left-symmetric algebra $A_{(3)}^1$.
Moreover it is a simple left-symmetric algebra, that is, it has not ideals
except itself and zero.

{\bf Proof}\quad The first half of conclusion follows directly from the
proof of Proposition 3.7 and the fact
that for every $c\ne 0$, $h(c,c)\ne 0$ since $h$ is the ordinary scalar
product.
The simplicity of the algebra is proved in [Bu2]. \hfill $\Box$

{\bf Remark 4}\quad The simple left-symmetric algebra $A_3^{(1)}$ is
firstly constructed in [Bu2]. In certain sense, our
re-construction gives it an interesting (geometric) interpretation.

Due to Corollary 2.7, we have

{\bf Corollary 4.2}\quad The scalar product appearing in equation (1.3) is
invariant under the sense of equation (2.10).

{\bf Corollary 4.3}\quad The (generalized) Burgers equation (1.4) is just
the following equation
$$u_t^1=u^1_{xx}+4u^1u_x^1+2\sum_{j=2}^{n}u^ju^j_x;\;\;u^k_t=u_{xx}^k+2u^1u^k_x-u^1u^1u^k-u^ku^ku^k,k=2,\cdots,
n.\eqno
(4.1)$$

{\bf Proof}\quad Let $C_{ij}^k$ be the structure constants. Hence equation
(1.4) gives
$$u^i_t=u^i_{xx}+2\sum_{j,k=1}^n C_{jk}^i u^ju^k_x+\sum_{k,j,l,m=1}^n
(C_{ml}^iC_{kj}^l-C_{kj}^lC_{lm}^i) u^ku^ju^m.$$
For the left-symmetric algebra $A_{(3)}^1$, the non-zero structure constants
are
$C_{11}^1=2, C_{jj}^1=1, C_{1j}^j=1,j=2,\cdots,n$. Hence equation (4.1)
follows.\hfill $\Box$

Besides the simple left-symmetric algebra $A_{(3)}^1$, there are some other
algebras appearing in Theorem 2.5 and Corollary 2.2 satisfying certain
additional (interesting) conditions, which play important roles in the study
of
left-symmetric algebras.

{\bf Definition} \quad Let $A$ be left-symmetric algebra.

(1) If for every $x\in A$, $R_x$ is nilpotent, then $A$ is said to be
transitive or complete. The transitivity corresponds to the completeness of
the affine manifolds in geometry ([Ki],[P]).

(2) If for every $x\in A$, $L_x$ is an interior derivation of the
sub-adjacent Lie algebra of $A$, then $A$ is said to be an interior
derivation
algebra. Such a structure corresponds to a flat left-invariant connection
adapted to the interior automorphism structure of a Lie group ([P]).

(3) If for every $x,y\in A$, $R_xR_y=R_yR_x$, then $A$ is said to be a
Novikov algebra. It was introduced in connection with the Poisson brackets
of hydrodynamic type and Hamiltonian operators in the formal variational
calculus ([GD],[BN]).

  (4) If for every $x,y,z\in A$, the associator $(x,y,z)$ is
right-symmetric, that is, $(x,y,z)=(x,z,y)$, then $A$ is said to be
bi-symmetric.
It is just the assosymmetric ring in the study of near-associative algebras
([Kl],[BM2]).

By direct computation, we have

{\bf Proposition 4.4}\quad Let $A$ be a left-symmetric constructed from
Theorem 2.5 and Corollary 2.2.

(1) $A$ is associative if and only if $A$ is isomorphic to one of the
following algebras: the associative algebras given in Corollary 2.2;
$A_{(1)}^{(k)}$; $A_{(2)}$; $A_{(7)}^1$.

(2) $A$ is transitive if and only if $A$ is trivial or $A$ is isomorphic to
$A_{(1)}^{(k)}$ or $A_{(4)}^0$.

(3) Besides the commutative cases, $A$ is an interior derivation algebra if
and only if $A$ is isomorphic to $A_{(4)}^0$. Moreover, $A_{(4)}^0$ is the
unique left-symmetric interior derivation algebra on the Lie algebra given
by equation (2.4). (cf. [P])

(4) Besides the commutative cases, $A$ is a Novikov algebra if and only if
$A$ is isomorphic to one of the following algebras: the associative algebra
in the case (2) of Corollary 2.2; $A_{(4)}^0$; $A_{(6)}$;
$A_{(7)}^\alpha$; $A_{(3)}^{(k),3}$.

(5) Besides the associative cases, $A$ is bi-symmetric if and only if $A$
isomorphic to $A_{(4)}^1$ or $A_{(6)}$.

{\bf Corollary 4.5}\quad Let $A$ be a left-symmetric constructed from
Theorem 2.5 and Corollary 2.2. Then $A$ with dimension $n>2$ is associative
(or transitive, or bi-symmetric, or a interior derivation algebra, or a
Novikov algebra) if and only if $A$ has such an additional structure when
its dimension $n=2$. Hence, the construction in this paper
can be regarded as generalization (not extension!) of certain 2-dimensional
left-symmetric algebras.

At the end of this paper, we give an application of the results in this
paper to integrable systems. Recall that a linear transformation $R$ on a
Lie algebra ${\cal G}$ is called a classical $r$-matrix if $R$ satisfies
$$[R(x), R(y)]=R([R(x),y]+[x,R(y)]),\;\;\forall x,y\in {\cal G}.\eqno
(4.2)$$
It corresponds to a solution of classical Yang-Baxter equation ([Ku4],[GS]).
Moreover, if $R$ satisfies the above equation, then
$$x*y=[R(x),y],\;\;\forall x,y\in {\cal G},\eqno (4.3)$$
defines a left-symmetric algebra on ${\cal G}$. Two classical $r$-matrices
are said to be equivalent if their corresponding left-symmetric algebras are
isomorphic.

{\bf Corollary 4.6}\quad For the Lie algebra $A$ given by equation
(2.4), there is only one non-zero classical $r$-matrix under the
sense of equivalence such that $A$ is the sub-adjacent Lie algebra
of left-symmetric algebra given by equation (4.3) , which $R$ is
given by
$$R(e_1)=e_1,\;\;R(e_j)=0,\;j=2,\cdots, n.\eqno (4.4)$$
The corresponding left-symmetric algebra given by equation (4.3) is
isomorphic to $A_{(4)}^0$.

{\bf Proof}\quad Let $R$ satisfy equation (4.2). Hence by equation (4.3), we
know that for every $x\in A$, $L_x={\rm ad}R(x)$, where ${\rm ad}$ is the
adjoint operator of Lie algebra. Hence $L_x$ is an interior derivation of
the Lie algebra $A$. Thus the left-symmetric algebra defined by equation
(4.3) is an interior derivation algebra. Therefore the conclusion follows
from (3) in Proposition 4.4. \hfill $\Box$

\newpage

\begin{center}{\bf ACKNOWLEDGMENTS}\end{center}

The author thanks Professor I.M. Gel'fand, Professor B.A. Kupershmidt and
Professor P. Etingof for useful
suggestion and great encouragement. The author also thanks Professor
J.
Lepowsky,
Professor Y.Z. Huang and Professor H.S. Li for the hospitality extended to
him during his stay at Rutgers, The State University of New Jersey and for
valuable discussions.
This work was supported in part by S.S. Chern Foundation for Mathematical
Research, the National Natural Science
Foundation of China(10201015), Mathematics Tianyuan Foundation(10126003),
K.C. Wong Education Foundation.

\begin{center}{\bf Reference}\end{center}

\baselineskip=20pt

\begin{description}

\item[[BM1]] C.M. Bai, D.J. Meng, The classification of left-symmetric
algebra in dimension 2, (in Chinese), Chinese Science
Bulletin 23 (1996) 2207.
\item[[BM2]] C.M. Bai, D.J. Meng, The structure of bi-symmetric algebras and
their sub-adjacent Lie algebras, Comm. in Algebra 28 (2000) 2717-2734.
\item[[BM3]] C.M. Bai, D.J. Meng, On the realization of transitive Novikov
algebras, J. Phys. A: Math. Gen. 34 (2001) 3363-3372.
\item[[BM4]] C.M. Bai, D.J. Meng, The realizations of non-transitive Novikov
algebras, J. Phys. A: Math. Gen. 34 (2001) 6435-6442.
\item[[BM5]] C.M. Bai, D.J. Meng, A Lie algebraic approach to Novikov
algebras, J. Geo. Phys. 45 (2003) 218-230.
\item[[BN]] A.A. Balinskii, S.P. Novikov, Poisson brackets of hydrodynamic
type, Frobenius algebras and Lie algebras, Soviet Math. Dokl. 32 (1985)
228-231.
\item[[Bu1]] D. Burde, Affine structures on nilmanifolds, Int. J. Math. 7
(1996) 599-616.
\item[[Bu2]] D. Burde, Simple left-symmetric algebras with solvable Lie
algebra, Manuscipta Math. 95 (1998) 397-411.
\item[[ES]] P. Etingof, A. Soloviev, Quantization of geometric classical
$r$-matrix, Math. Res. Lett. 6 (1999) 223-228.
\item[[ESS]] P. Etingof, T. Schedler, A. Soloviev, Set-theoretical solutions
to the quantum Yang-Baxter equations, Duke Math. J. 100 (1999)169-209.
\item[[GD]] I.M. Gel'fand,  I. Ya. Dorfman, Hamiltonian operators and
algebraic structures related to them, Functional Anal. Appl. 13 (1979)
248-262.
\item[[GS]] I.Z. Golubschik, V.V. Sokolov, Generalized operator Yang-Baxter
equations, integrable ODES and nonassociative algebras, J. Nonlinear Math.
Phys., 7 (2000) 184-197.
\item[[Ki]] H. Kim, Complete left-invariant affine structures on nilpotent
Lie groups, J. Diff. Geo. 24 (1986) 373-394.
\item[[Kl]] E. Kleinfeld, Assosymmetric rings, Proc. Amer. Math. Soc. 8
(1957) 983-986.
\item[[Ku1]] B.A. Kuperschmidt, Non-abelian phase spaces, J. Phys. A: Math.
Gen. 27 (1994)2801-2810.
\item[[Ku2]] B.A. Kuperschmidt, Quantum differential
forms,  J. Nonlinear
Math. Phys., 5 (1998) 245-288.
\item[[Ku3]] B.A. Kuperschmidt, On the nature of the Virasoro algebra, J.
Nonlinear Math. Phys., 6 (1999) 222-245.
\item[[Ku4]] B.A. Kuperschmidt, What a classical $r$-matrix really is, J.
Nonlinear Math. Phys., 6 (1999) 448-488.
\item[[P]] A.M. Perea, Flat left-invariant connections adapted to the
automorphism structure of a Lie group, J. Diff. Geo. 16
(1981) 445-474.
\item[[S]] S.I. Svinolupov, On the analogues of the Burgers equation, Phys.
Lett. A 135 (1989) 32-36.
\item[[SS]] S.I. Svinolupov, V.V. Sokolov, Vector-matrix generalizations of
classical integrable equations, Theoret. and Math. Phys. 100
(1994) 959-962.
\item[[V]] E.B. Vinberg, Convex homogeneous cones, Transl. of Moscow Math.
Soc. No. 12 (1963) 340-403
\end{description}

\end{document}